\documentclass[12pt]{article}
\usepackage{amsfonts, amsmath, amssymb, latexsym}

\numberwithin{equation}{section}
\newcommand{\ad}{\hbox{\rm ad\,}}
\newcommand{\Ann}{\hbox{\rm Ann\,}}
\newcommand{\Hom}{\hbox{\rm Hom\,}}
\newcommand{\Der}{\hbox {\rm Der}}
\newcommand{\eqdef}{\buildrel{\hbox{\small def}} \over =}

\newcommand{\qed}{\quad $\square$}

\newtheorem{Thm}{Theorem}[section]
\newtheorem{Pro}[Thm] {Proposition}

\begin{document}

\title{Non-degenerate graded Lie algebras\\
with a degenerate transitive subalgebra}
\author{T. B. Gregory \\ M. I. Kuznetsov\footnote{The second author gratefully acknowledges 
support from the Russian Foundation of Basic Research Grant
\#{}05-01-00580.   \newline}}
\date{}
 \maketitle

\bigskip

\begin{abstract} The property of degeneration of modular graded Lie algebras,
first investigated by B. Weisfeiler, is analyzed. Transitive irreducible graded
Lie algebras $L=\sum _{i\in \mathbb Z}L_i,$ over algebraically closed fields 
of characteristic $p>2,$ with classical reductive component $L_0$ are considered. 
We show that if a non-degenerate Lie algebra $L$ contains a transitive degenerate 
subalgebra $L'$ such that $\dim L'_1>1,$ then $L$ is an infinite-dimensional Lie algebra.
\end{abstract}

\setcounter{section}{-1}
\begin{section}{Introduction} \end{section}

One of the most important steps the program of classifying the simple finite-dimensional Lie algebras
of characteristic $p>0$, a program developed by A.I. Kostrikin and I.R. Shafarevich [10], is the investigation of non-contractible filtrations of simple Lie algebras. Let $L_0$ be a maximal subalgebra in a simple Lie 
algebra $L$ such that the nilradical of the adjoint representation of $L_0$ on $L$ is nontrivial, let $L_{-1}$ be an $L_0$-submodule of $L$ such that $L_0\subset L_{-1},$ and suppose that $L_{-1}/L_0$ is an irreducible
$L_0$-submodule. The non-contractible filtration of $L$ corresponding to the pair $(L_{-1}, L_0)$
is constructed by induction:
$$L_{-i}=[L_{-1}, L_{-i+1}]+L_{-i+1}, L_i=\{ l\in L_{i-1} | [l, L_{-1}]\subset L_{i-1}\}, i>0.$$

It follows from the definition of a non-contractible filtration that the associated graded Lie algebra
${\mathfrak g}=\oplus _{i=-q}^r{\mathfrak g}_i$ possesses the following properties:

\bigskip

\begin{enumerate}

\item[{\rm (1)}] (irreducibility) ${\mathfrak g}_{-1}$ is an irreducible  ${\mathfrak g}_0$-module;

\item[{\rm (2)}] (transitivity) for any $j\geqq 0$ if $x\in {\mathfrak g}_j$ and
$[x,{\mathfrak g}_{-1}]=(0)$, then $x=0$;

\item[{\rm (3)}] ${\mathfrak g}_{-i}=[{\mathfrak g}_{-i+1},{\mathfrak g}_{-1}]$ for any $i>1;$ 

\end{enumerate}

A graded Lie algebra ${\mathfrak g}=\oplus _{i\in \mathbb Z}{\mathfrak g}_i$ satisfying conditions
(1) - (3) is called a transitive irreducible Lie algebra. A graded subalgebra
${\mathfrak g}'\subset {\mathfrak g}$ is called {\it transitive} if 
${\mathfrak g}'_{-1}={\mathfrak g}_{-1}$.  
In [14] and [15],  B. Weisfeiler investigated the properties of finite-dimensional transitive
irreducible graded Lie algebras $\mathfrak g$ over algebraically closed fields of characteristic 
$p>0$ and obtained results of fundamental importance for the theory of modular graded Lie algebras.  
He showed ([15]) that the subalgebra ${\mathfrak g}^-=\sum_{i<0}{\mathfrak g}_i$ contains a unique maximal ideal $M=M({\mathfrak g})$ of $\mathfrak g$ called the Weisfeiler radical of $\mathfrak g$, such that the factor algebra $\mathfrak g/M$ is a semisimple Lie algebra with a unique minimal ideal $I$.
 The centroid of the ideal $I$ is a truncated polynomial algebra. 

Lie algebras of the form ${\mathfrak g}$
fall into two classes: non-degenerate and degenerate Lie algebras. In the non-degenerate case,
the centroid of the ideal $I$ has zero degree, and the grading of ${\mathfrak g}/M$ is determined by
the grading of a simple Lie algebra, namely the core of the differentially simple ideal $I$. In the degenerate case, the grading in $\mathfrak g$ is determined by a nontrivial grading of the centroid; in addition,
 ${\mathfrak g}_2=0$ and $ [[{\mathfrak g}_{-1},{\mathfrak g}_1], {\mathfrak g}_1]=0$.
Weisfeiler showed ([15], Proposition 3.2.1) that the last property is a criterion for degeneration in 
finite-dimensional transitive irreducible Lie algebras $\mathfrak g$. 

When
 investigating the transitive irreducible graded Lie algebras generated by the local part
 ${\mathfrak g}_{-1}+{\mathfrak g}_0+{\mathfrak g}_1,$  it is not {\it a priori} clear whether the algebra 
${\mathfrak g}$ is finite dimensional or not. Therefore, we assume Weisfeiler's criterion of degeneration of finite-dimensional Lie algebras to be the definition of a degenerate (resp. non-degenerate) transitive
irreducible $\mathbb Z$-graded Lie algebra ${\mathfrak g}$.
Properties (1) -- (3) are asymmetric with respect to ${\mathfrak g}_{-1}$ and ${\mathfrak g}_1$.
It is therefore natural to consider the subalgebra  ${\mathfrak g}'$ in $\mathfrak g$ generated by the local part 
${\mathfrak g}_{-1}+{\mathfrak g}_0+{\mathfrak g}'_1$ where ${\mathfrak g}'_1$ is an irreducible  ${\mathfrak g}_0$-submodule of ${\mathfrak g}_1.$ (Note that ${\mathfrak g}'$ satisfies the above definition of a transitive subalgebra.) However, even if $\mathfrak g$ is non-degenerate,   ${\mathfrak g}'$ might nonetheless be degenerate. 
We investigate this problem  for the case in which ${\mathfrak g}_0$ is a classical reductive Lie algebra. This case
is of particular interest for the classification theory of simple modular Lie algebras. The Recognition Theorem describing such Lie algebras satisfying the additional condition of transitivity with respect to ${\mathfrak g}_1$, was obtained by V. Kac  [8] for $p>5$ under assumption that 
 ${\mathfrak g}_{-1}$ is a restricted ${\mathfrak g}_0$-module. G. M. Benkart, T.B. Gregory and A. A. Premet [2] extended the Recognition Theorem for the case $p>3$.  When $p=3,$ there exist Lie algebras which are unlike any in characteristics $p>3$ yet 
satisfy the conditions of the Recognition Theorem (see [3], [5], [12]).  
In the present paper, a theorem (Theorem 2.1) is proved that is of great importance for 
the classification of graded Lie algebras with a classical reductive component  
${\mathfrak g}_0$ and nonrestricted ${\mathfrak g}_0$-module ${\mathfrak g}_{-1}$ when $p>2$.

\begin{Thm} 
Let $L=\oplus _{i\in {\mathbb Z}} L_i$ be a non-degenerate transitive irreducible graded Lie algebra
with classical reductive component  $L_0$ over an algebraically closed field ${\mathbf F}$ of characteristic $p>2$. 
If $L$ contains a degenerate transitive subalgebra $L'$ such that  $\dim L'_1>1$, then  $L$ is an infinite-dimensional Lie algebra. 
\end{Thm}

By transitivity, the representation of $L_0'$ on $L_1$ is restricted when and only when the representation of $L_0'$ on $L_{-1}$ is restricted.  Since no non-restricted representation of $L_0'$ can have dimension 1, we have the following corollary.

\bigskip

\noindent {\bf Corollary 0.2}
Let $L$ be as in the above theorem, and suppose that the representation of $L_0'$ on $L_{-1}$ is not restricted. Then  $L$ is an infinite-dimensional Lie algebra.

\bigskip

It should be pointed out that both the case in which $\dim L'_1=1$ and the case of even characteristic are currently being investigated elsewhere.

The proof of Theorem 0.1 is given in Section 2. Section 1 contains needed definitions, notations, and results obtained by B. Weisfeiler.

\begin{section}{Preliminaries} \end{section}

 Recall that {\it a classical Lie algebra} over a field ${\mathbf F}$ of
characteristic $p>0$ can be obtained from a $\mathbb Z$-form (the ``Chevalley basis'')
of a complex simple Lie algebra by reducing the scalars modulo $p$ and
 extending them to ${\mathbf F}$. This process may result in a Lie algebra with a non-zero
center; such a Lie algebra is still referred to as ``classical'',
as is the quotient of such a Lie algebra by its center. For
example, the Lie algebras $\mathfrak{sl}(pk)$ and
$\mathfrak{psl}(pk)$ are both considered to be classical Lie
algebras. It could also happen that a classical Lie algebra has a noncentral
ideal, as does the Lie algebra G$_2$ if $p=3$. The Lie algebras $\mathfrak{gl}(pk)$
and $\mathfrak{pgl}(pk)$ are also considered classical.

{\it A classical reductive Lie algebra\/} $\mathfrak{g}$ is the
sum of commuting ideals $\mathfrak{g}_j$ which are classical Lie
algebras, and an at-most-one-dimensional center
$\mathfrak{z}(\mathfrak{g})$:

\begin{equation}\mathfrak{g} = \mathfrak{g}_1 + \cdots + \mathfrak{g}_k + \mathfrak{z}(\mathfrak{g}) \nonumber \end{equation}

Let ${\mathcal O}_n={\mathbf F}[x_1,\ldots ,x_n]/(x_1^p,\ldots ,x_n^p)$ be a truncated 
polynomial algebra. Denote by $W_n$ the Lie algebra of vector fields $\Der {\mathcal O}_n.$
The following theorem of B. Weisfeiler [15] plays a fundamental
r\^ole in the study of graded Lie algebras.

\bigskip

\begin{Thm} {\bf (Weisfeiler's Theorem)}   
Let $L=L_{-q}\oplus \cdots \oplus L_{-1} \oplus
 L_0 \oplus L_1 \oplus \cdots \oplus L_r$  be a transitive irreducible graded Lie algebra
over an algebraically closed field ${\mathbf F}$ of characteristic $p>0,$ and let $M(L)$ be the largest
ideal of $L$ contained in   $L_{-q}\oplus \cdots \oplus L_{-1}.$ Then

(i) $L/M(L)$ is semisimple and contains a unique minimal ideal 
$I=S\otimes {\mathcal O}_n$, where $S$ is a simple Lie algebra and
$n$ is a non-negative integer. The ideal $I$ is a graded ideal,  and $I_i=(L/M(L))_i$ for any $i<0.$

(ii) (Degenerate case) If $I_1=(0),$ then there exists an integer $k$, $1 \leqq k  \leqq n$ such that 
 the grading of ${\mathcal O}_n,$ induced by the grading of
$I,$ is given by setting $\deg (x_i)=-1$ for  
$1\leqq i\leqq k,$ and $\deg (x_i)=0$ for $k  < i\leqq n.$
Moreover,   $I_i=S\otimes {\mathcal O}_{n,i}$ for any $i$; 
$L_2=(0)$; $I_0=[L_{-1},L_1]$; $L_1\subseteq \{D\in 1\otimes W_n
\vert \deg (D)=1\}$; $[[L_{-1},L_1],L_1]=0$; and
 $$L_0\subset \Der(S\otimes {\mathcal O}_{n- k}) +
1\otimes {\mathcal O}_{n-k}\otimes W_{k,0}, 
W_{k,0}\cong {\mathfrak {gl}}({\mathcal O}_{ k,-1}).$$

(iii) (Non-degenerate case) If  $I_1\ne (0),$ then $S$ is a graded Lie algebra, and 
$I_i=S_i\otimes {\mathcal O}_n$ for any $i$. Moreover, 
$(0) \ne [L_{-1},L_1] \subseteq I_0.$  
If $\Der S=\oplus (\Der S)_i$ is the grading of the Lie algebra $\Der S$
induced by the grading of $S,$ then 
\begin{equation} L_0\subset (\Der S)_0\otimes {\mathcal O}_n + 
1\otimes W_n\nonumber\end{equation}
\noindent  и для $i>0$ $G_i\subset (\Der S)_i\otimes {\mathcal O}_n.$
\end{Thm}

A transitive graded Lie algebra $L=\oplus _{i\in \mathbb Z}L_i$ is called {\it degenerate\/} if
$[[L_{-1},$ $L_1],$ $L_1]$ $=0.$ The following proposition proved by B. Weisfeiler ([15], Proposition 3.2.1) motivates this definition.

\begin{Pro} 
Let $L$ be a finite-dimensional transitive irreducible graded Lie algebra, and let $V$ be a $L_0$-submodule of $L_1$. Suppose that  $L_0$ is not faithful on $V.$ Then

(i) $[V,V]=0,$ and $[[L_{-1},V],V]=0.$

(ii) Let $L'$ be the subalgebra of $L$ generated by $L_{-1}+L_0+V.$ Then $L'/M(L')$ satisfies the  
conditions and conclusions of the degenerate case of Theorem 1.1.
\end{Pro}

Set  $L^- = \oplus_{i<0}L_i$ and $ L^+ = \oplus_{i>0}L_i.$
Let $L=L_{-q}+\ldots +L_{-1}+L_0+L_1+\ldots $ be a $q$-graded Lie algebra and set
$G=\oplus _{i\in \mathbb Z} L_{iq}.$ Then $G$ is a graded subalgebra of $L$ which can be considered to be a 1-graded Lie algebra if we set $G_i=L_{iq}$. Denote by $\mathcal T$ the largest ideal of $G$ contained in $G_0+G^+$. The factor algebra $G/\mathcal T$  is denoted by $B(L_{-q})$. Evidently, $B(L_{-q})$ is a transitive 1-graded Lie algebra.  
The design of $B(L_{-q})$ is a particular case of the construction given in [1].

\begin{section}{The proof of the Theorem} \end{section} 

To simplify notation, we formulate Theorem 0.1 another way.

\bigskip

\begin{Thm}  Let $L=\oplus _{i\in {\mathbb Z}} L_i$ be a non-degenerate transitive irreducible graded
Lie algebra over an algebraically closed field ${\mathbf F}$ of characteristic $p>2$ with classical reductive component  $L_0.$ If  $[[L_{-1}, V], V] = 0$ for some $L_0$-submodule
$V \subset L_1$ such that $\dim V > 1,$ then $\dim L = \infty $. \end{Thm}
The {\bf Proof} of the Theorem consists of several steps.
We suppose that $\dim L < \infty $ and obtain a contradiction.
In $Aut \, L,$ denote by $T$ the one-dimensional torus which defines the
$\mathbb Z$-grading of $L$. All algebras constructed below will be $T$-invariant, and thus inherit the ${\mathbb Z}$-grading. Note that the symbols $L^{\dagger}$, $L_i^{\dagger},$ etc., used below, retain their meanings only within a particular proof segment (i.e., {\it (a)}, {\it (b)}, {\it (c)}, etc.)

\bigskip

\noindent {\it (a)  If $V_1$ is an $L_0$-submodule in $L_1$, $V\subset V_1$ and
$[[L_{-1}, V_1], V_1] = 0,$ then $V_1 = V$\/.}

Let $L^{\dagger}$ be the subalgebra of $L$ generated by the local Lie algebra $L_{-1} + L_0 + V_1$. By Weisfeiler's Theorem, 
$L^{\dagger}$ is a degenerate Lie algebra and $[L_{-1}, V_1]$ is a minimal ideal of $L_0$. 
It follows from Theorem 1.1 (ii)  that $S = [L_{-1}, V_1]$ is a simple ideal of $L_0$ and that $L_{-1}= S \otimes {\mathcal O}_{n,-1},$ 
where ${\mathcal O}_{n,-1} = \, <x_1,\ldots ,  x_n>\,$ is an irreducible $L_0/S$-module. 
Since $V_1 = \, <\partial _1,\ldots,\partial _n>\, = \, <x_1,\ldots ,x_n>^*,$ it follows that $V_1$ is an irreducible $L_0/S$-module 
and, therefore, an irreducible
$L_0$-module. Hence, $V_1 = V,$ so (a) is proved.

\bigskip

\noindent {\it (b) Here the problem is reduced to the case where
$L_1/V$ is an irreducible $L_0$-module, $L^+$ is generated by $L_1$ and
$M(L)=0$. The subalgebra ${\mathcal L}_0$ generated by $L_{-1} + L_0 + V$
is described\/.}  

Since by assumption $L$ is  non-degenerate, i.e., $[[L_{-1},L_1],L_1]\neq 0,$ it follows that $L_1 \neq V.$ Let $\mathfrak{V}$ be an
$L_0$-submodule of $L_1$ such that $V \subset \mathfrak{V}$ and $0 \neq \mathfrak{V}/V$
is an irreducible $L_0$-module. It follows from (a) that
$[[L_{-1},\mathfrak{V}],\mathfrak{V}] \neq 0.$ We will consider the subalgebra of $L$ generated by the local part $L_{-1} + L_0 + \mathfrak{V}$. Since we are assuming that $L$ is finite dimensional, this subalgebra is finite dimensional, also, and
it satisfies the hypotheses of Theorem 2.1. Therefore, we will prove Theorem 2.1 by replacing the original $L$ (if necessary) by this subalgebra, 
which we will henceforth refer to as $L$. Now, out of this new $L$ we can factor its Weisfeiler radical $M(L)$. 
Thus, we will consider a non-degenerate finite-dimensional Lie algebra $L$ satisfying the conditions 
of Theorem 2.1 such that $L$ is generated by its local part, $L_1/V$ is an irreducible $L_0$-module, and $M(L)=0$.

Denote by $\mathcal L_0$ the subalgebra of $L$ generated by $L_{-1} + L_0 + V$.
As noted in the proof of part (a) above, it follows from Weisfeiler's Theorem and the assumption that the null component is classical reductive that $S = [L_{-1}, V]$ is a simple
ideal of $L_0$. According to (a), $V$ is an irreducible $L_0$-module.
Evidently, $M({\mathcal L}_0)$ 
 is both a maximal $V$-invariant ideal of
$L^{-}$ and $T$-invariant. By Weisfeiler's Theorem, $[V, V] = 0$ and

$${\mathcal L}_0/M({\mathcal L}_0) = L^{\dagger}_{-|\delta |} + L^{\dagger}_{-|\delta | + 1} + \ldots + L^{\dagger}_{-1} + L_0 + V,$$

\noindent where

$$L^{\dagger}_{-1} = L_{-1}, \, \hbox{and} \, L^{\dagger}_{-i} = S\otimes {\mathcal O}_{n, -i}, i > 0; $$ 

\noindent here $\delta$ is the n-tuple $(p-1, \ldots, p-1) $ and ${\mathcal O}_n 
= \oplus _{i \leqq 0}{\mathcal O}_{n,i}$ is the grading 
of ${\mathcal O}_n={\mathbf F}[x_1, \ldots ,x_n]/(x_1^p, \ldots ,x_n^p)$ opposite to the standard grading. 
  Since ${\mathcal O}_{n,-1} =  \, <x_1,\ldots ,x_n> \, \cong V^*,$ it follows that $n = \dim V > 1.$ 
Futhermore, $L_{-1} = S\otimes {\mathcal O}_{n,-1},$ and $V =id_S\otimes <\partial_1,\ldots ,\partial_n>\,$ 
$\subset \Hom_{\mathbf F}(L_{-1}, L_{0}).$ We will identify $V$ with
$\, <\partial _1,\ldots ,\partial _n>\,$. By Weisfeiler's Theorem, $S \subset
L_0 \subset (\Der S)\otimes 1 + 1\otimes W_{n,0}\subset \Hom_{\mathbf F}(L_{-1},L_{-1})$,
where $W_{n,0} = \, <x_i\partial _j, i,j=1,\ldots , n>\,$ and 1 denotes the identity
map. Let $\pi :(\Der S)\otimes 1 + 1\otimes W_{n,0} \longrightarrow W_{n,0}$ be
the projection along the ideal $(\Der S)\otimes 1.$ Since $L_{-1}$ is an
irreducible $L_0$-module, it follows that $<x_1,\ldots , x_n>\,$ is an irreducible $\pi (L_0)$-module. 

\bigskip

\noindent {\it (c) Here it is shown that ${\mathcal L}_0$ is a maximal subalgebra of $L$ and that $L/{\mathcal L}_0$ 
contains a unique nontrivial irreducible ${\mathcal L}_0$-submodule 
$\overline{{\mathcal L}_{-1}}={\mathcal L}_{-1}/{\mathcal L}_0$. The non-contractible filtration of $L$ corresponding to
${\mathcal L}_{-1}$ and ${\mathcal L}_0$ is constructed, as is its associated graded Lie algebra
${\mathfrak g}$\/.}

It follows from Weisfeiler's Theorem that for $i > 1,$ $\{l\in L_i\vert  (\ad
L_{-1})^{i-1}l\subset V \}=0,$ since otherwise there would exist a degenerate Lie algebra
$L_{-q}+\ldots +L_{-1}+L_0 + V_1 +V_2+\ldots $ in which $V_1=V$ and $V_2=\{ l\in L_2\vert (\ad L_{-1})l\subset V \} \neq 0.$ 
Since $L$ is a transitive graded Lie algebra and $L_1/V$ is an irreducible $L_0$-module,
it follows that any nontrivial ${\mathcal L}_0$-submodule of $L/{\mathcal L}_0$ contains
the unique irreducible ${\mathcal L}_0$-submodule of $L/{\mathcal L}_0$, namely, the intersection
of all ${\mathcal L}_0$-submodules of $L/{\mathcal L}_0$ which contain $(L_1 + {\mathcal L}_0)/{\mathcal L}_0$. 
Since $L^+$ is generated by $L_1$, it follows that ${\mathcal L}_0$ is a maximal subalgebra of $L$.

Denote by ${\mathcal L}_{-1}$ the ${\mathcal L}_0$-submodule of $L$ such that
${\mathcal L}_0\subset {\mathcal L}_{-1}$ and ${\mathcal L}_{-1}/{\mathcal
L}_0$ is the unique irreducible ${\mathcal L}_0$-submodule of $L/{\mathcal
L}_0$ described in the preceding paragraph. Note that $L, {\mathcal L}_0, L/{\mathcal L}_0$ and, therefore,
${\mathcal L}_{-1}$ are invariant under the torus $T$, so ${\mathcal L}_{-1}$ is 
a ${\mathbb Z}$-graded subspace of $L$. Let

$$L={\mathcal L}_{-s}\supset {\mathcal L}_{-s+1}\supset \ldots \supset {\mathcal
L}_{-1}\supset {\mathcal L}_0\supset {\mathcal L}_1\supset \ldots $$

\noindent be the non-contractible  filtration of $L$ corresponding to the pair $({\mathcal
L}_{-1}, {\mathcal L}_0)$ where ${\mathcal L}_{-i} = {\mathcal L}_{-1}^i + {\mathcal
L}_{-i+1}$ for $i>1$ and ${\mathcal L}_j=\{ l\in {\mathcal L}_{j-1} \vert
[{\mathcal L}_{-1}, l]\subset {\mathcal L}_{j-1} \}$ for $j>0.$ 
Obviously, the filtration $\{ {\mathcal L}_i\}$ is invariant with respect to the
torus $T$. 

Let ${\mathfrak g} = gr L = \oplus _{i=-s}^r{\mathfrak g}_i$ be the graded Lie
algebra associated with the filtration $\{ {\mathcal L}_i\}$. Since the
filtration is $T$-invariant, it follows that $T$ acts on ${\mathfrak g}$ by
homogeneous automorphisms. Thus, ${\mathfrak g}_i=
\oplus_{j\in {\mathbb Z}}{\mathfrak g}_{ij}$ and ${\mathfrak g}=\oplus_{i,j\in
{\mathbb Z}}{\mathfrak g}_{ij}$ is a bigrading of ${\mathfrak g}$. The second
index of ${\mathfrak g}_{ij}$ refers to its elements' weight with respect to the torus $T$.
Evidently, ${\mathcal L}_1$ is a maximal
$\ad_L$-nilpotent ideal of ${\mathcal L}_0$; therefore, $M({\mathcal L}_0)\subset
{\mathcal L}_1$, where $M({\mathcal L}_0)$ is the Weisfeiler radical of the
degenerate graded Lie algebra ${\mathcal L}_0$. By Weisfeiler's Theorem
${\mathcal L}_0/M({\mathcal L}_0)$ is a semisimple Lie algebra, so that
$M({\mathcal L}_0) = {\mathcal L}_1,$ and

$${\mathfrak g}_0={\mathcal L}_0/{\mathcal L}_1=
{\mathcal L}_0/M({\mathcal L}_0)= S\otimes {\mathcal O}_n + L_0 + V.$$

\noindent The ${\mathbb Z}$-gradings in ${\mathfrak g}_0$ and in 
${\mathcal O}_n$ corresponding to the torus $T$ coincide with the grading corresponding to the
degenerate graded Lie algebra ${\mathcal L}_0/M({\mathcal L}_0)$ in Weisfeiler's
Theorem.

\bigskip

\noindent {\it (d) Here we show that ${\mathfrak g}_{-1}$ has the following ${\mathbb
Z}$-grading corresponding to the torus $T$:

$$ {\mathfrak g}_{-1}={\mathfrak g}_{-1,1}+\ldots +{\mathfrak
g}_{-1,|\delta |+1},  |\delta |=n(p-1),$$

\noindent where ${\mathfrak g}_{-1,i} = (\ad V)^{i-1}(L_1/V)$. Moreover, ${\mathfrak
g}_{-1}$ is a graded ${\mathcal O}_n$-module
 where the grading in 
${\mathcal O}_n$ is induced by the torus $T$, $\deg x_i = -1, i = 1,\ldots,n$, the bracket operation
of $S\otimes {\mathcal O}_n\subset {\mathfrak g}_0$ with ${\mathfrak g}_{-1}$ is
${\mathcal O}_n$-bilinear, and ${\mathfrak g}_{-1, |\delta | + 1}$ is both an irreducible 
${\mathfrak g}_{0,0}$-module and a nontrivial $S$-module. Here ${\mathfrak
g}_{0,0}$ is the 0-term of the ${\mathbb Z}$-grading of ${\mathfrak g}_0$
corresponding to the torus $T$, and $S=S\otimes 1\subset L_0={\mathfrak g}_{0,0}$\/.} 
 
Denote by ${\mathfrak h}$ the subalgebra of ${\mathfrak g}_0$ equal to $S\otimes
{\mathcal O}_n +L_0.$ Then ${\mathfrak g}_0=V\oplus {\mathfrak h}.$ Note that
$\ad_L V$ consists of nilpotent elements and that $(\ad _{{\mathfrak g}_0}l)^p=0$ for any
$l\in V.$ Since ${\mathfrak g}_{-1}$ is an irreducible ${\mathfrak g}_0$-module, 
it follows that $(\ad_{{\mathfrak g}_{-1}}l)^p=0$ for any $l\in V.$ 
Let $A=U({\mathfrak g}_0)/<l^p, l\in V>\,$ be the quotient of the universal
enveloping algebra of ${\mathfrak g}_0$ by the ideal generated by the set 
$\{l^p, l\in V \},$ and set $B=U({\mathfrak h}).$ Since $L_1/V$ is an irreducible $L_0-$module, 
it is therefore also an irreducible ${\mathfrak h}$-module; consequently, ${\mathfrak g}_{-1}$ is covered by the induced
${\mathfrak g}_0$-module $A\otimes _B(L_1+{\mathcal L}_0)/{\mathcal L}_0={\mathbf F}[V]\otimes _{\mathbf F}L_1/V$ where
${\mathbf F}[V]={\mathbf F}[\partial _1, \ldots ,\partial _n]/(\partial _1^p, \ldots ,\partial
_n^p)$ is a truncated polynomial algebra. Here $\{ \partial _1,\ldots ,\partial
_n\}$ is a basis of $V.$ It follows that
$${\mathfrak g}_{-1}=\sum_{i=0}^{ | \delta | }(\ad V)^i(L_1/V)=L_1/V\oplus \ad
V(L_1) \oplus \ldots \oplus (\ad V)^{ | \delta | }(L_1).$$
In particular, the ${\mathbb Z}$-grading of ${\mathfrak g}_{-1}$ is as follows
$${\mathfrak g}_{-1}=\oplus _{i=1}^{ |\delta |+1} {\mathfrak g}_{-1,i}, 
{\mathfrak g}_{-1,i}=(\ad V)^{i-1}(L_1/V).$$

We now show that ${\overline L}_1=L_1/V$ is a nontrivial $S$-module. Since
${\overline L}_1$ is an irreducible $L_0$-module and $S$ is an ideal of $L_0,$
it follows that $[S, {\overline L}_1]$ is equal to ${\overline L}_1$ or zero. In the latter
case $\ad S(L_1)\subset V,$ and since $[S, V]=0$, it follows that $(\ad_{L_1}s)^2=0$
for any $s\in S,$ so $\ad_{L_1}S$ is nilpotent by Engel's Theorem. But $S$ is simple, 
so, by a version of Schur's Lemma, either $\ad_{L_1}S \cong S$ or $\ad_{L_1}S = 0;$ 
therefore, being nilpotent,  $\ad_{L_1}S=0.$ Hence,
$\ad l_1:L_{-1}=S\otimes {\mathcal O}_{n,-1}\longrightarrow L_0$ is a nontrivial
morphism of $S$-modules for any $0\neq l_1\in L_1;$ thus, $[L_{-1}, L_1]$ is an 
isotypical $S$-submodule of $L_0$ of type $S.$ Therefore, $[L_{-1}, L_1]=S$
and $[[L_{-1}, L_1], L_1]=[S, L_1]=0.$  We have obtained a contradiction, since $L$ is
a nondegenerate graded Lie algebra. Thus, $\ad S({\overline L}_1)={\overline
L}_1.$ Since $[S,V]=0$, it follows that for any $i$ 
\begin{equation} \label {C} [S, {\mathfrak g}_{-1,i}]=[S, (\ad V)^{i-1}({\overline L}_1)]=(\ad V)^{i-1}[S,
{\overline L}_1]={\mathfrak g}_{-1,i}.
\end{equation}

Since $S\otimes {\mathcal O}_n$ is a minimal ideal of ${\mathfrak g}_0,$ 
it follows that $(\Ann_{{\mathfrak g}_0}{\mathfrak g}_{-1})\cap (S\otimes {\mathcal O}_n)=0$.
Therefore, since $deg S\otimes x^{\delta }=-|\delta |$ and $[S\otimes
x^{\delta }, {\mathfrak g}_{-1}]\neq 0$, it follows, in view of the 
${\mathbb Z}$-grading in ${\mathfrak g}_{-1},$ that 

\begin{equation}\label {A} [S\otimes x^{\delta }, {\mathfrak g}_{-1, |\delta |+1}]\neq 0. \end{equation} 

\noindent Hence, ${\mathfrak g}_{-1, |\delta
|+1}\neq 0.$ Since ${\mathfrak g}_{-1, |\delta |+1}=(ad V)^{|\delta
|}{\overline L}_1=\ad(\partial _1^{p-1}\cdot \ldots \cdot \partial
_n^{p-1}){\overline L}_1$, it follows that the $L_0$-module ${\mathfrak g}_{-1, |\delta |+1}$
is isomorphic to the twisted $L_0$-module ${\overline L}_1\otimes {\mathbf F}_{\sigma }$
where ${\mathbf F}_{\sigma }$ is a one-dimensional $L_0$-module corresponding to
a morphism $\sigma :L_0\longrightarrow {\mathbf F}.$ Thus, ${\mathfrak g}_{-1, |\delta
|+1}$ is an irreducible $L_0$-module.

For any $B$-module $U$ denote by $c(U)$ the truncated coinduced ${\mathfrak g}_
0$-module $coind_B^AU=\Hom_B(A, U),$ where $A$ and $B$ are as above.  Let $
{\tilde {\mathfrak g}}_{-1}=\sum_{i=1}^{|\delta |}{\mathfrak g}_{-1,i}.$
Obviously, ${\tilde {\mathfrak g}}_{-1}$ is an ${\mathfrak h}$-submodule of the 
${\mathfrak g}_0$-module ${\mathfrak g}_{-1}$ and ${\mathfrak g}_{-1}/ {\tilde
{\mathfrak g}}_{-1}\cong {\mathfrak g}_{-1, |\delta |+1}.$
According to the theory of truncated coinduced modules ([11], [13]), 
there exists  a nonzero morphism of ${\mathfrak g}_0$-modules
from ${\mathfrak g}_{-1}$ to $c({\mathfrak g}_{-1, |\delta |+1})$ which is
injective since ${\mathfrak g}_{-1}$ is an irreducible ${\mathfrak g}_0$-module.
Furthermore, ${\mathcal O}_n=c({\mathbf F})$ where ${\mathbf F}$ is a trivial ${\mathfrak h}$-module, and
$S\otimes {\mathcal O}_n=c(S\otimes {\mathcal O}_n/S\otimes {\mathfrak n})=
c(S)$ where $\mathfrak n$ is the maximal ideal of ${\mathcal O}_n$. In addition, any module $c(U)$ is a free ${\mathcal O}_n$-module. The
bracket operation $[ , ]$ in ${\mathfrak g}$, $[ , ]:S\times {\mathfrak g}_{-1,
|\delta |+1}\longrightarrow {\mathfrak g}_{-1, |\delta |+1}$ induces the
${\mathcal O}_n$-bilinear mapping

$$\mu :c(S)\times c({\mathfrak g}_{-1, |\delta |+1})\longrightarrow
c({\mathfrak g}_{-1, |\delta |+1}),$$

\noindent $\mu (\phi ,\psi ) =[ , ]\circ \phi \otimes \psi \circ \delta,$ where $\delta $ is the
coproduct in $A.$
As in Proposition 2.2 of [7], it may be shown that the
bracket operation in the Lie algebra ${\mathfrak g}, [ , ]:S\otimes {\mathcal
O}_n\times {\mathfrak g}_{-1}\longrightarrow {\mathfrak g}_{-1}$ coincides
with the restriction of the mapping $\mu $. Note that a basis of the space
${\mathfrak g}_{-1, |\delta |+1}$ is a basis of the free ${\mathcal
O}_n$-module $c({\mathfrak g}_{-1, |\delta |+1}).$ 

Since the bracket operation of $S$ $\otimes$ ${\mathcal O}_n$ with ${\mathfrak g}_{-1}$
is  ${\mathcal O}_n$-bilinear and $[S,$ ${\mathfrak g}_{-1, |\delta
|+1}]$ $={\mathfrak g}_{-1, |\delta |+1},$ it follows that

\begin{eqnarray}c({\mathfrak g}_{-1, |\delta |+1}) &=& {\mathcal O}_n{\mathfrak g}_{-1, |\delta
|+1}\nonumber \\
&=& {\mathcal O}_n[S, {\mathfrak g}_{-1, |\delta |+1}]=[S\otimes {\mathcal
O}_n, {\mathfrak g}_{-1, |\delta |+1}]\nonumber \\
&\subseteq& {\mathfrak g}_{-1}\subseteq
c({\mathfrak g}_{-1, |\delta |+1}).\nonumber\end{eqnarray}

Thus, ${\mathfrak g}_{-1}=c({\mathfrak g}_{-1, |\delta |+1})={\mathfrak g}_{-1,
|\delta |+1}\otimes {\mathcal O}_n, {\mathfrak g}_{-1,i}={\mathfrak g}_{-1,
|\delta |+1}\otimes {\mathcal O}_{n,i-|\delta|-1}.$
\bigskip
 
{\bf Remark.} It may be inferred from the results of R. Block [4] that ${\mathfrak
g}_{-1}$ is a free ${\mathcal O}_n$-module and that the bracket operation of ${\mathfrak
g}_{-1}$ with $S\otimes {\mathcal O}_n$ is ${\mathcal O}_n$-bilinear.

\bigskip

\noindent {\it (e) Here it is shown that ${\mathfrak g  }_1\neq 0\/.$}

According to (c), ${\mathcal L}_1=M({\mathcal L}_0),$ where ${\mathcal
L}_0=L_{-q}+\ldots +L_{-1}+L_0+V.$ Suppose that $M({\mathcal L}_0)=0$. Then

$${\mathcal L}_0={\mathfrak g}_0=S\otimes {\mathcal O}_n + L_0 + V, 
$$

\noindent $q=
|\delta |,$ $L_{-|\delta |}$ $= S\otimes x^{\delta },$ and $L_{-|\delta |-1}=0.$ It follows from (d) that the 
subspaces ${\mathfrak g}_{-1,i},$ $i=2, \ldots ,|\delta |+1$ may be identified 
with subspaces $L^{\dagger}_i$ $\subset L_i;$ moreover,
$$[L_{-1}, L^{\dagger}_{|\delta |+1}]=[S\otimes {\mathcal O}_{n,-1}, {\mathfrak g}_{-1, |\delta |+1}]=
{\mathcal O}_{n,-1} {\mathfrak g}_{-1,|\delta |+1}={\mathfrak g}_{-1,|\delta |}=L^{\dagger}_{|\delta |}.         $$

We now show that $[[L_{-|\delta |}, L^{\dagger}_{|\delta |}], L^{\dagger}_{|\delta |}]\neq 0$ and that
$[[L_{-|\delta |}, L^{\dagger}_{|\delta |}], L_{-|\delta |}]\neq 0$. Bearing \eqref{A} in mind, let 
$U\eqdef [L_{-|\delta |}, {L^{\dagger}}_{|\delta |+1}]\subset L_1.$ Since the bracket operation of 
${\mathfrak g}_{-1}$ with $S\otimes {\mathcal O}_n\subset {\mathfrak g}_0$ is ${\mathcal
O}_n$-bilinear, it follows that $U/U\cap V=L_1/V={\mathfrak g}_{-1,1}$ (see (d)). By the
assumptions on $L,$ $V$ and $L_1/V$ are irreducible $L_0$-modules, so either
$V\subset U,$ or $L_1 = V \oplus U$ is a direct sum of $L_0$-modules. Furthermore,

\begin{equation} \label {B}[L_{-|\delta |}, L^{\dagger}_{|\delta |}]=[L_{-|\delta |}, [L_{-1}, L^{\dagger}_{|\delta |+1}]]=
[L_{-1}, [L_{-|\delta |}, {L^{\dagger}}_{|\delta |+1}]]=[L_{-1}, U].
\end{equation}

\noindent We next show that $S\subset [L_{-1}, U].$ Suppose $[L_{-1},U]\cap S=0.$ Since
$[L_{-1}, U]$ and $S$ are ideals of $L_0,$ and $S$ is simple, it follows that $[[L_{-1}, U],
S]=0.$ Since $[L_{-1}, V]=S$, it follows that $V\cap U=0$; therefore, $L_1=U\oplus V$.

Consider the subalgebra $L^{{\dagger}\dagger}$ in $L$ generated by $L_{-1} + L_0 + U$,
 
$$L^{{\dagger}\dagger}=L_{-|\delta |}+ L_{-|\delta |+1}+ \ldots +L_{-1}+L_0+U+ \ldots .$$

\noindent Now $M(L^{{\dagger}\dagger})$ is an ideal of $S\otimes {\mathcal O}_n=L_{-|\delta |}+ \ldots
+L_{-1}+S$;  therefore, $M(L^{{\dagger}\dagger})=S\otimes J$ where $J$ is an ideal of ${\mathcal
O}_n$. Inasmuch as $M(L^{{\dagger}\dagger})\subset L_{-|\delta |}+\ldots +L_{-2}$, we have $J\neq
{\mathfrak m}$ where ${\mathfrak m}$  is the maximal ideal of ${\mathcal O}_n$. Evidently,

\begin{eqnarray} \overline {L^{{\dagger}\dagger}} &=& L^{{\dagger}\dagger}/M(L^{{\dagger}\dagger})=\overline {L}_{-k}+ \ldots +\overline {L}_{-2}
+L_{-1}+L_0 +U +\ldots, \nonumber \\
\overline {L}_i  &=& S\otimes ({\mathcal O}_n/J)_i, i=-2,
\ldots , -k.\nonumber\end{eqnarray}

Suppose $\overline {L^{{\dagger}\dagger}}$ is a nondegenerate Lie algebra. Since $L_0$ is a
reductive classical Lie algebra, it follows from Weisfeiler's Theorem that $\overline {L^{{\dagger}\dagger}}$
contains a simple graded Lie algebra 

\begin{eqnarray} S^{{\dagger}{\dagger}} &=& \oplus _iS^{{\dagger}{\dagger}}_i, S^{{\dagger}{\dagger}}\subset \overline
{L^{{\dagger}\dagger}}\subset \Der S^{{\dagger}{\dagger}},\nonumber \\
 \overline {L^{{\dagger}\dagger}}_i &=& S^{{\dagger}{\dagger}}_i, i=-1, \ldots ,-k, S^{{\dagger}{\dagger}}_j\subset
L^{{\dagger}\dagger}_j\subset \Der_j S^{{\dagger}{\dagger}}, j=0,1,\ldots .\nonumber\end{eqnarray}

\noindent Note that $S^{{\dagger}{\dagger}}_{-k}$ is an irreducible $S^{{\dagger}{\dagger}}_0$-module and
 
$$S^{{\dagger}{\dagger}}_0=[S^{{\dagger}{\dagger}}_{-1}, S^{{\dagger}{\dagger}}_1]=[S^{{\dagger}{\dagger}}_{-1},U]=[L_{-1},U].$$

\noindent Let $\rho :L_0\longrightarrow {\mathfrak gl}(S^{{\dagger}{\dagger}}_{-k})$ be the restriction of
the adjoint representation. By Schur's Lemma, $\dim C_{{\mathfrak
gl}(S^{{\dagger}{\dagger}}_{-k})}(\rho (S^{{\dagger}{\dagger}}_0))=1$.
On the other hand $S^{{\dagger}{\dagger}}_{-k}=\overline {L^{{\dagger}\dagger}}_{-k}=S\otimes ({\mathcal
O}_n/J)_{-k}.$
Inasmuch as $S\cap S^{{\dagger}{\dagger}}_0=S\cap [L_{-1}, U]=(0)$, we have $S\cong \rho (S)\subset
C_{{\mathfrak gl}(S^{{\dagger}{\dagger}}_{-k})}(\rho (S^{{\dagger}{\dagger}}_0))$ and we obtain a contradiction.

If $\overline {L^{{\dagger}\dagger}}$ is a degenerate Lie algebra, then by Weisfeiler's Theorem,
$\overline {L^{{\dagger}\dagger}}$ contains a differentially simple ideal
$S^{{\dagger}{\dagger}}\otimes {\mathcal O}_m,$ ${\mathcal O}_m$ $={\mathbf F}[y_1,\ldots ,y_m]/(y_1^p,\ldots
,y_m^p),$ $ \overline {L^{{\dagger}\dagger}}_{-i}$ $=S^{{\dagger}{\dagger}}\otimes 
{\mathcal O}_{m,-i},$ $i=-1,\ldots ,-k $ where ${\mathcal O}_m$ $=\oplus {\mathcal
O}_{m,-i}$ is the grading of ${\mathcal O}_m$ opposite to the standard grading
and $S^{{\dagger}{\dagger}}$ $=[L_{-1}, U].$ In particular, $\overline {L^{{\dagger}\dagger}}_{-k}$ $=S^{{\dagger}{\dagger}}\otimes y^{\delta
'},$ where $\delta ' =(p-1,\ldots ,p-1)$ is an $m$-tuple. Therefore, since we determined above that $[[L_{-1}, U],
S]=0,$ it follows that $\overline {L^{{\dagger}\dagger}}_{-k}$
is a trivial $S$-module. On the other hand,  $\overline {L^{{\dagger}\dagger}}_{-k}=S\otimes
({\mathcal O}_n/J)_{-k}$ is an isotypical $S$-module of the type $S$, and we
obtain a contradiction. Thus, $S\subset [L_{-1}, U].$

Since $L_{-|\delta |}=S\otimes x^{\delta },$ it follows (from \eqref{B}) that

$$[[L_{-|\delta |}, L^{\dagger}_{|\delta |}], L_{-|\delta |}]=[[L_{-1}, U], 
L_{-|\delta| }] \supseteq 
[S, L_{-|\delta |}]=L_{-|\delta |}\neq 0.$$

\noindent According to (d), $L_1/V$ is a nontrivial $S$-module and $L^{\dagger}_{|\delta |}$ $=
[[\ldots[L_1/V, \underbrace{V], \ldots ],V}_{|\delta |-1}],$ whence $[S, L^{\dagger}_{|\delta
|}]=L^{\dagger}_{|\delta |}\neq 0.$ Therefore,
 
$$[[L_{-|\delta |},L^{\dagger}_{|\delta |}], L^{\dagger}_{|\delta |}]=[[L_{-1}, U], 
L^{\dagger}_{|\delta|}] \supseteq 
[S, L^{\dagger}_{|\delta |}]=L^{\dagger}_{|\delta |}.$$

\noindent Consider the subalgebra $B(L_{-|\delta |})$ (see Section 1). Since
 
$$ [[L_{-|\delta |}, L^{\dagger}_{|\delta |}], L_{-|\delta |}] = L_{-|\delta |} {\rm \, and \, } 
     [[L_{-|\delta |}, L^{\dagger}_{|\delta |}], L^{\dagger}_{|\delta |}] = L^{\dagger}_{|\delta |},$$

\noindent it follows that the Lie algebra $B(L_{-|\delta |})$ is a one-graded transitive irreducible Lie
algebra,

$$ B(L_{-|\delta |})=B_{-1}+B_0+B_1+\ldots, B_1\neq 0.$$

\noindent Since $B_{-1}=L_{-|\delta |}=S\otimes x^{\delta }$, it follows that $S\subset B_0\subset \Der
S +{\mathfrak z}(B_0), B_{-1}\cong S.$ According to the theorem of Kostrikin-Ostrik
([9]) such an algebra does not exist when $p>2$.

Thus, $M({\mathcal L}_0)\neq 0$ and, therefore, ${\mathfrak g}_1\neq 0$.

\bigskip

\noindent {\it (f)} Let $L=\ldots \supset {\mathcal L}_{-1}\supset {\mathcal
L}_0\supset 
{\mathcal L}_1\supset \ldots $ be the non-contractible filtration of $L$
corresponding to the maximal subalgebra ${\mathcal L}_0=L_{-q}+\ldots +L_{-1}
+L_0 +V$ and the ${\mathcal L}_0$-module ${\mathcal L}_{-1}$,

\begin{equation}
\label{filtr}
{\mathcal L}_{-1}={\mathcal L}_0+L_1+\sum _{i=1}^{|\delta |}(ad V)^iL_1 
\end{equation}

\noindent (see (c), (e)).  Here the following statements are proved:

\bigskip

\noindent  {\it 1) let $l\in {\mathcal L}_1$ and suppose that $[l, V]=0;$ then $l\in {\mathcal
L}_2$ if and only if $[l, L_1]\subset {\mathcal L}_1$;

\bigskip

\noindent 2) let ${\mathfrak g}_1=\oplus _{i=k}^l {\mathfrak g}_{1, -i},
{\mathfrak g}_{1,-k}\neq 0$ be the ${\mathbb Z}$-grading of ${\mathfrak g}_1$
corresponding to the torus $T;$ then $k=min \{i| M({\mathcal L}_0)\cap L_{-i}\neq
0 \},$ and we have $2\leqq k\leqq |\delta |+1$;
$${\mathcal L}_1\subset L_{-q}+\ldots +L_{-k}, {\mathcal L}_2\subset
L_{-q}+\ldots +
L_{-k-1},$$
${\mathfrak g}_{1,-k}\cong M({\mathcal L}_0)\cap L_{-k}$;

\bigskip

\noindent 3) $[{\mathfrak g}_1, {\mathfrak g}_{-1}]=S\otimes 
{\mathcal O}_n\subset {\mathfrak g}_0$\/.}

\bigskip

According to (c), ${\mathcal L}_1=M({\mathcal L}_0)\subset L_{-q}+\ldots +L_{-2}$.
Evidently, if $l\in {\mathcal L}_2,$ then $[l, L_1]\subset[l, {\mathcal
L}_{-1}]\subset 
{\mathcal L}_1.$ Let $l\in {\mathcal L}_1,$ and suppose that $[l, V]=0$. According to
\eqref{filtr},
 
$$ [l, {\mathcal L}_{-1}]=[l, {\mathcal L}_0] +[l, L_1]+
\sum _{i=1}^{|\delta |}(ad V)^i([l, L_1]).$$

\noindent Therefore, if $[l, L_1]\subset {\mathcal L}_1,$ then $[l, {\mathcal L}_{-1}]$ is
contained in ${\mathcal L}_1$; that is, $l\in {\mathcal L}_2,$ and 1) is proved.

To prove 2), set $s=min \{i| M({\mathcal L}_0)\cap L_{-i}\neq 0\}$. Then $k\geq
s\geq 2.$ Let
$L^{\dagger}_{-s}=M({\mathcal L}_0)\cap L_{-s}$ and $L^{{\dagger}\dagger}_{-s}
=\{ l\in L^{\dagger}_{-s}| [l, L_1]=0\}$. Obviously, $L^{{\dagger}\dagger}_{-s}$ is an $L_0$-submodule of $L^{\dagger}_{-s}$. Since $L^+$ is
generated by $L_1$, it follows that the ideal $J$ of $L$ generated by $L^{{\dagger}\dagger}_{-s}$ is contained in $L_{<}$. By the assumptions on $L,$
$M(L)=0$; therefore, $J=0$ and $L^{{\dagger}\dagger}_{-s}=0$. Hence, $0\neq [l,L_1]\subset
L_{-s+1}$ for any $0\neq l\in L^{\dagger}_{-s}\subset {\mathcal L}_1$. By the definition
of $s,$ $[l, V]= M({\mathcal L}_0) \cap L_{-s+1} = 0$ for any $l\in L^{\dagger}_{-s}$; therefore, according to 1), $l\notin
{\mathcal L}_2$ for any $0\neq l\in L^{\dagger}_{-s}$. Thus, $k=s$ and ${\mathfrak
g}_{1,-k}\cong L^{\dagger}_{-k}=M({\mathcal L}_0)\cap L_{-k}.$  Since the depth of the 
${\mathbb Z}$-grading of the semisimple Lie algebra ${\mathcal L}_0/{\mathcal
L}_1={\mathcal L}_0/ M({\mathcal L}_0)$ is equal to $|\delta |$, it follows that 
$L_{-|\delta |-1}\subset M({\mathcal L}_0).$  Therefore, $k\leqq |\delta |+1.$

\noindent 3) Since $V\subset {\mathcal L}_0$ and ${\mathcal L}_1$ is a
${\mathcal L}_0$-module, we have from \eqref{filtr} that 

$$ [{\mathcal L}_1, {\mathcal L}_{-1}]\subset {\mathcal L}_1 +[{\mathcal L}_1,
L_1]+
\sum _{i>0}(adV)^i[{\mathcal L}_1,L_1].$$

\noindent As  ${\mathcal L}_1\subset L_{-q}+\ldots +L_{-2}$, we have that $[{\mathcal L}_1,
L_1]\subset
L_{-q}+\ldots +L_{-1}$ and $(ad V)^i[{\mathcal L}_1, L_1]\subset L_{-q} +\ldots
+L_{-1} +S,$  

$$([{\mathcal L}_1, {\mathcal L}_{-1}] +{\mathcal L}_1)/{\mathcal L}_1\subset
S\otimes {\mathcal O}_n\subset {\mathfrak g}_0={\mathcal L}_0/{\mathcal L}_1;$$

\noindent that is, $[{\mathfrak g}_1, {\mathfrak g}_{-1}]\subset S\otimes {\mathcal O}_n.$
Since $S\otimes {\mathcal O}_n$ is a minimal ideal of ${\mathfrak g}_0$ and
$[{\mathfrak g}_{1, -k}, {\mathfrak g}_{-1,1}]\cong [L^{\dagger}_{-k}, L_1]\neq 0,$ it follows that
$[{\mathfrak g}_1, {\mathfrak g}_{-1}] = S\otimes {\mathcal O}_n.$

\bigskip

\noindent {\it (g) Let $k$ be as in (f). We show that $[[{\mathfrak g}_{-1}, 
{\mathfrak g}_1],
  {\mathfrak g}_1]\neq 0$ and that if $k=|\delta |+1$, then 
${\mathfrak g}_{1, -k}$ is a nontrivial $S$-module, and ${\mathfrak g}_{1,-|\delta
|-1}\cong 
L_{-|\delta |-1}$\/.}

Note that if $[S, {\mathfrak g}_{1, -k}]\neq 0$, then according to (f), 3) 
$[[{\mathfrak g}_{-1}, {\mathfrak g}_1], {\mathfrak g}_1]\neq 0.$ Suppose that 
$[S, {\mathfrak g}_{1,-k}]=0.$

Let $k\leqq |\delta |.$ As in (f), we set $L^{\dagger}_{-k}={\mathcal L}_1\cap L_{-k}$ and
show that 
$L^{\dagger}_{-k-1}=[L_{-1}, L^{\dagger}_{-k}]\neq 0,  L^{\dagger}_{-k-1}\not\subset {\mathcal L}_2.$
Since 
$L^{\dagger}_{-k}\cong  {\mathfrak g}_{1,-k}, L_{-1}=S\otimes {\mathcal O}_{n, -1}\cong
{\mathfrak g}_{0, -1}\subset 
[{\mathfrak g}_{-1},{\mathfrak g}_1],$ it will follow that 
$[S\otimes {\mathcal O}_{n, -1}, {\mathfrak g}_{1, -k}]\neq 0;$ that is, 
$[[{\mathfrak g}_{-1}, {\mathfrak g}_1], {\mathfrak g}_1]\neq 0.$

Since $[V, L^{\dagger}_{-k}]$ $=0$, it follows that $[V,$ $L^{\dagger}_{-k-1}]$ $= [V,$ $[L_{-1},$ $L^{\dagger}_{-k}]]$ $= 
[[V,$ $L_{-1}],$ $L^{\dagger}_{-k}]$ $=[S,$ $L^{\dagger}_{-k}]$ $\cong [S,$ ${\mathfrak g}_{1,-k}]$ $=0.$
Furthermore,

$$[L_1, [L_{-1}, L^{\dagger}_{-k}]] \subseteq [[L_1, L_{-1}], L^{\dagger}_{-k}] + [L_{-1}, [L_1, L^{\dagger}_{-k}]].$$

As $[L_0, L^{\dagger}_{-k}]\subset L^{\dagger}_{-k}\subset {\mathcal L}_1$ and 
$0\neq [L_1, L^{\dagger}_{-k}]\subset L_{-k+1}=S\otimes {\mathcal O}_{n,-(k-1)}$, we have that

\begin{eqnarray} 0\neq  ([L_{-1}, [L_1, L^{\dagger}_{-k}]]+{\mathcal L}_1)/{\mathcal L}_1 &=& 
([L_{-1}, ~[L_1, ~L^{\dagger}_{-k}]]+L^{\dagger}_{-k})/~L^{\dagger}_{-k}\nonumber \\
&\subset& [S\otimes {\mathcal O}_{n,-1}, ~S\otimes {\mathcal O}_{n, -k+1}]\nonumber \\
&=& S\otimes {\mathcal O}_{n, -k}.\nonumber\end{eqnarray}

\noindent So $[[L_{-1}, L^{\dagger}_{-k}], L_1]\not\subset {\mathcal L}_1$ and by (f), 1)
$[L_{-1},L^{\dagger}_{-k}]=L^{\dagger}_{-k-1}\not\subset {\mathcal L}_2.$

Now let $k>|\delta |.$ Then $k=|\delta |+1, L^{\dagger}_{-k}=L_{-k}\cong {\mathfrak
g}_{1,-k}, 
L_{-i}= S\otimes {\mathcal O}_{n,-i}, i=1, \ldots , |\delta |,
L_{-k}=[L_{-1}, L_{-|\delta |}]= [S\otimes {\mathcal O}_{n,-1}, S\otimes
x^{\delta }].$
By assumption, $[S, L_{-k}]=0.$ For $x\in {\mathcal O}_{n, -1},$ the
bracket operation in $L$,
$[ , ]: S\otimes x\times S\otimes x^{\delta }\longrightarrow L_{-k}$ gives us a
mapping of the
$S$-module $S\otimes S$ into the trivial $S$-submodule $[S\otimes x, S\otimes
x^{\delta }].$
Since $S$ is a classical simple Lie algebra, $S\cong S^*$ as $S$-modules.
This fact follows, in particular, from Curtis's Theorem. (See [6] for $p>7,$
and [9] for any $p>0.$)
Therefore, the quotients  of  $S\otimes S$ which are trivial over $S$, are in
one-to-one correspondence
with the trivial $S$-submodules in $S^*\otimes S\cong (S\otimes S^*)^*=\Hom(S, S).$
By Schur's Lemma $\Hom (S, S)$ has the unique nonzero trivial $S$-module $\Hom_S(S,
S), \dim \Hom_S(S, S)=1.$  So, $[S\otimes x, S\otimes x^{\delta }]= \, <\hat{x}>\,$ is
a one-dimensional trivial $S$-submodule in $L_{-|\delta |-1}$ and

$$[s_1\otimes x, s_2\otimes x^{\delta }]=\varphi _x(s_1, s_2)\hat{x},$$

\bigskip

\noindent where $\varphi _x(s_1, s_2)$ is an invariant bilinear form on $S$. Since $\Phi
_0=\{ x\in
{\mathcal O}_{n, -1}| \varphi _x=0 \}$ is a $\pi (L_0)$-submodule of ${\mathcal
O}_{n, -1}$ and 
${\mathcal O}_{n, -1}$ is an irreducible $\pi (L_0)$-module, it follows that $\Phi _0=0.$
(See the description of $\pi (L_0)$ in (b).)

Let $\{ x_1, \ldots ,x_n \}$ be a basis of ${\mathcal O}_{n, -1}$. The above
discussion shows that the space of invariant bilinear forms on $S$ is
one-dimensional. So, the elements $\hat{x_i}$ may be chosen in such a way that
$[s_1\otimes x_i, s_2\otimes x^{\delta }]=
\phi (s_1, s_2) \hat{x_i}$ for a fixed nonzero invariant form $\phi $ on $S$.

Obviously, $L_{-|\delta |-1}$ is isomorphic to ${\mathcal O}_{n, -1}$ as a $\pi
([L_0, L_0])$-module, and is therefore an irreducible $L_0$-module. Let
$L^{\dagger}_{|\delta |+1}=
(\ad V)^{|\delta |}L_1$. According to (d), $L^{\dagger}_{|\delta |+1}$ is not only not zero but also
 isomorphic to the irreducible ${\mathfrak g}_{0,0}$-module 
${\mathfrak g}_{-1, |\delta |+1}$. Note that ${\mathfrak g}_{0,0}=L_0$. Since
$S$ is an ideal of $L_0$, it follows that $L^{\dagger}_{|\delta |+1}$ is an isotypical $S$-module.
Since $[L_{-|\delta |-1}, L_1]=
L_{-|\delta |}=S\otimes x^{\delta }$ and $[L_{-|\delta |-1}, V]=0,$ it follows that

\begin{eqnarray} [L_{-|\delta |-1}, L^{\dagger}_{|\delta |+1}] &=& 
[L_{-|\delta |-1}, (\ad V)^{|\delta |}L_1]\nonumber \\
&=& (\ad V)^{|\delta |}[L_{-|\delta |-1}, L_1]= 
(\ad V)^{|\delta |}S\otimes x^{\delta }= S.\nonumber\end{eqnarray}

\noindent Now, as $L_{-|\delta |-1}$ is an irreducible $L_0$-module which is trivial over $S$, the
mapping $\ad l, l\in L_{-|\delta |-1}$ is a nonzero $S$-morphism from
$L^{\dagger}_{|\delta |+1}$ to $S$. Therefore, 
$L^{\dagger}_{|\delta |+1}$ is an isotypical $S$-module of type $S$ and may be
represented as $S\otimes U$ where $U$ is an irreducible $L_0$-module which is trivial
over $S$; that is, $U$ is an irreducible $\pi (L_0)$-module. The bracket operation
in $L$ of elements 
$\hat{x}\in L_{-|\delta |-1}$ and $s\otimes u\in L^{\dagger}_{|\delta |+1}=S\otimes U$
may be written as
$[\hat{x}, s\otimes u]=\, < \hat{x}, u>s$ where $< , >$ is a nondegenerate pairing,
$< , >\,:L_{-|\delta |-1}\times L^{\dagger}_{|\delta |+1}\longrightarrow {\mathbf F}.$ 

Let $\alpha $ be a root of the classical simple Lie algebra $S,$ let $e_{\alpha },
e_{-\alpha }$ be root vectors, and let $h_{\alpha }=[e_{\alpha }, e_{-\alpha }],$ so that $sl(2)=
\, <e_{\alpha }, e_{-\alpha }, h_{\alpha }>\,$ is the corresponding three-dimensional simple
subalgebra of $S$. Evidently, $\varphi (e_{-\alpha }, e_{\alpha })\neq 0.$ For 
$e_{\alpha }\otimes x^{\delta }\in  L_{-|\delta |}$ and $e_{\alpha }\otimes u\in
S\otimes U=L^{\dagger}_{|\delta |+1},$  let $\overline{[e_{\alpha }\otimes x^{\delta
},e_{\alpha }\otimes u]}$ be the corresponding coclass in $L_1/V={\mathfrak
g}_{-1,1}.$ Since the bracket of 
$S\otimes {\mathcal O}_n\subset {\mathfrak g}_0$ with $ {\mathfrak g}_{-1}$ is 
${\mathcal O}_n$-bilinear (see (c)),
$$\overline{[e_{\alpha }\otimes x^{\delta }, e_{\alpha }\otimes
u]}=\overline{[e_{\alpha }, e_{\alpha }\otimes u]}x^{\delta
}=\overline{[e_{\alpha }, e_{\alpha }]\otimes u}x^{\delta }=0.$$
Consequently, in $L$ we have $[e_{\alpha }\otimes x^{\delta }, e_{\alpha }\otimes u]\in V\subset
L_1.$
Inasmuch as $[S, V]=0$, we have 
$[h_{\alpha },[e_{\alpha }\otimes x^{\delta }, e_{\alpha }\otimes u]]=0$.
However,

$$[h_{\alpha },[e_{\alpha }\otimes x^{\delta }, e_{\alpha }\otimes u]]=
2\alpha (h_{\alpha })[e_{\alpha }\otimes x^{\delta }, e_{\alpha }\otimes u].$$

\noindent Thus, $[e_{\alpha }\otimes x^{\delta }, e_{\alpha }\otimes u]=0.$

Let $<\hat{x},u>\,\neq 0$.  Then in $L$ we have

\begin{eqnarray} 0 &=& [e_{-\alpha }\otimes x, [e_{\alpha }\otimes x^{\delta }, e_{\alpha }\otimes
u]]\nonumber \\
 &=& \varphi (e_{-\alpha }, e_{\alpha })[\hat{x}, e_{\alpha }\otimes u]\nonumber\\
&&\qquad +
[e_{\alpha }\otimes x^{\delta }, [e_{-\alpha }, e_{\alpha }]\otimes ux]\nonumber \\
&=& \varphi (e_{-\alpha },e_{\alpha })\, <\hat{x}, u>\,e_{\alpha }- 
[e_{\alpha }\otimes x^{\delta }, h_{\alpha }\otimes ux].\nonumber\end{eqnarray}

\noindent Therefore,

$$0\neq [e_{\alpha }\otimes x^{\delta }, h_{\alpha }\otimes ux]=\varphi
(e_{-\alpha }, e_{\alpha })<\hat{x}, u>\,e_{\alpha }\in S\subset L_0.$$

\noindent Thus, $S\subset [L_{-|\delta |}, L^{\dagger}_{|\delta |}]$ where $L^{\dagger}_{|\delta |}={\mathfrak g}_{-1, |\delta |}$ (see (d)). Note that $L_{-|\delta |}=S\otimes x^{\delta }$. It now follows from \eqref{C}
that
 
$$[[L_{-|\delta |}, L^{\dagger}_{|\delta |}], L^{\dagger}_{\pm |\delta |}]\supset [S, L^{\dagger}_{\pm
|\delta |}]= L_{\pm |\delta |},$$

\noindent where $L^{\dagger}_{-|\delta |}=L_{-|\delta |}.$

Consider the subalgebra ${\mathfrak h}$ in ${\mathfrak g}$ generated by the local
part 
${\mathfrak h}_{-1}+{\mathfrak h}_0+{\mathfrak h}_1,$ where 
${\mathfrak h}_{-1}={\mathfrak g}_{0,-|\delta |}=S\otimes x^{\delta }\cong
L_{-|\delta |},
{\mathfrak h}_0={\mathfrak g}_{0,0}=L_0,$ and ${\mathfrak h}_1=
{\mathfrak g}_{-1, |\delta |}\cong L^{\dagger}_{|\delta |}=(S\otimes U)\otimes {\mathcal
O}_{n,-1}.$
Evidently, ${\mathfrak h}$ is a one-graded Lie algebra with respect to the 
${\mathbb Z}$-grading corresponding to the torus $T$.  The transitive one-graded
Lie algebra $B=B_{{\mathfrak h}}({\mathfrak h}_{-1})$ satisfies the conditions
of the
 Kostrikin-Ostrik Theorem [9]. Obviously, $B=B_{-1}+B_0+B_1+\ldots,$ where
$B_1\neq 0,
S\subset B_0\subset \Der S+{\mathfrak z}(B_0),$ and $B_{-1}\cong S.$ It follows from the
Kostrikin-Ostrik Theorem that such an algebra does not exist. We have arrived at a
contradiction. Therefore,
${\mathfrak g}_{1,-k}\cong L_{-k}$ is a nontrivial $S$-module in the case 
$k=|\delta |+1$. Therefore, $[[ {\mathfrak g}_{-1},{\mathfrak g}_1],{\mathfrak
g}_1]\neq 0.$

\bigskip

\noindent {\it (h) In this section of the proof, it is proved that there exists a simple graded
Lie algebra $\hat{S}=
\oplus _{i\in {\mathbb Z}}\hat{S}_i$ such that 
${\mathfrak g}_{-1}=\hat{S}_{-1}\otimes {\mathcal O}_n,
\hat{S}_0\otimes {\mathcal O}_n\subset {\mathfrak g}_0\subset
(\Der _0\hat{S})\otimes {\mathcal O}_n+1\otimes W_n,
\hat{S}_1\otimes {\mathcal O}_n\subset {\mathfrak g}_1\subset (\Der
_1\hat{S})\otimes {\mathcal O}_n,$ and $\hat{S}_0=S$. The bracket operation in
${\mathfrak g}, 
[ , ]: {\mathfrak g}_{-1}\times {\mathfrak g}_1\longrightarrow [{\mathfrak
g}_{-1},
{\mathfrak g}_1]=S\otimes {\mathcal O}_n$ is the restriction of the ${\mathcal
O}_n$-bilinear
bracket operation of the Lie algebra $(\Der \hat{S})\otimes {\mathcal O}_n , 
[{\mathfrak g}_{-1}, {\mathfrak g}_1]=\hat{S}_0\otimes {\mathcal O}_n$\/.}

Consider the subalgebra $\tilde{{\mathfrak g}}$ of  ${\mathfrak g}$ generated by
${\mathfrak g}_{-1}, {\mathfrak g}_0, {\mathfrak g}_1$. The torus $T$ acts by
automorphisms on  the subalgebra  $\tilde{{\mathfrak g}};$ therefore,
$\tilde{{\mathfrak g}}$ has a ${\mathbb Z}$-grading corresponding to the torus
$T$. Let $\hat{{\mathfrak g}}=
\tilde{{\mathfrak g}}/M(\tilde{{\mathfrak g}})$ where $M(\tilde{{\mathfrak g}})$
is the Weisfeiler radical of $\tilde{{\mathfrak g}}$. According to (g), 
$[[ {\mathfrak g}_{-1},{\mathfrak g}_1],{\mathfrak g}_1]\neq 0$. So
$\hat{{\mathfrak g}}$ is a nondegenerate graded Lie algebra. By Weisfeiler's
Theorem, $\hat{{\mathfrak g}}$ contains
a minimal ideal $A(\hat{{\mathfrak g}})=\hat{S}\otimes{\mathcal O}_m$ where $\hat{S}$
is a graded simple Lie algebra, $\hat{S}=\oplus \hat{S}_i,  \hat{{\mathfrak
g}}_i=\hat{S}_i\otimes {\mathcal O}_m $ for $i<0,$ $\hat{{\mathfrak
g}}_{-1}={\mathfrak g}_{-1},
 \hat{{\mathfrak g}}_0={\mathfrak g}_0, \hat{{\mathfrak g}}_1={\mathfrak g}_1,
\hat{S}\otimes {\mathcal O}_m\subset {\mathfrak g}_0\subset
(\Der_0\hat{S})\otimes {\mathcal O}_m +
1\otimes W_m=\Der_0(\hat{S}\otimes {\mathcal O}_m),$ and
$\hat{S}_1\otimes {\mathcal O}_m\subset {\mathfrak g}_1\subset
(\Der_1\hat{S})\otimes {\mathcal O}_m.$
Since $A(\hat{{\mathfrak g}})$ is the unique minimal ideal of 
$\hat{{\mathfrak g}}$, it is invariant with respect to the torus $T$.  
 
According to (f), 3), 
$[{\mathfrak g}_{-1}, {\mathfrak g}_1]=S\otimes {\mathcal O}_n.$
On the other hand, since $[\hat{S}_{-1},\hat{S}_1]=\hat{S}_0$ in the graded simple
Lie algebra $\hat{S}$, it follows that $[{\mathfrak g}_{-1}, {\mathfrak
g}_1]=\hat{S}_0\otimes {\mathcal O}_m;$ the bracket operation

$$ [{\mathfrak g}_{-1}, {\mathfrak g}_1]=\hat{S}_0\otimes {\mathcal O}_m = 
S\otimes {\mathcal O}_n$$

\noindent is
${\mathcal O}_m$-bilinear; and ${\mathcal O}_m$ is naturally contained in
the centroid ${\mathfrak C}$ of the Lie algebra $S\otimes {\mathcal O}_n$. Inasmuch as ${\mathfrak C}$  is isomorphic to ${\mathcal O}_n$, we have ${\mathcal O}_m \subset {\mathcal O}_n$.
Since ${\mathcal O}_m$ is isomorphic to the centroid of the $T$-invariant ideal 
$\hat{S}\otimes {\mathcal O}_m$, it follows that ${\mathcal O}_m$ is a graded subalgebra in
${\mathcal O}_n$ with respect to the ${\mathbb Z}$-grading corresponding to $T$.
  
 Let ${\mathfrak m}$ be the maximal ideal of ${\mathcal O}_m$, and let
${\mathfrak n}$ be the maximal ideal of ${\mathcal O}_n$.  We have the following
series of natural isomorphisms:

$$\hat{S}_0\cong 
(\hat{S}_0\otimes _{\mathbf F} {\mathcal O}_m)\otimes _{{\mathcal O}_m}
{\mathcal O}_m/~{\mathfrak m}{\mathcal O}_m =
(S\otimes _{\mathbf F}{\mathcal O}_n)\otimes _{{\mathcal O}_m}{\mathcal O}_m/{\mathfrak m}
{\mathcal O}_m$$

$$= S\otimes _{\mathbf F}({\mathcal O}_n\otimes_{{\mathcal O}_m} {\mathcal O}_m/{\mathfrak m}{\mathcal
O}_m)\cong
S\otimes _{\mathbf F}{\mathcal O}_n/{\mathfrak m}{\mathcal O}_n.$$

\noindent Therefore, $\hat{S}_0\otimes _{\mathbf F}{\mathcal O}_m\cong S\otimes _{\mathbf F}B\otimes _{\mathbf F} 
{\mathcal O}_m$ where $B={\mathcal O}_n/{\mathfrak m}{\mathcal O}_n$. Hence,
${\mathcal O}_n\cong B\otimes _{\mathbf F} {\mathcal O}_m$ and so $B\cong {\mathcal O}_l,
{\mathcal O}_n\cong {\mathcal O}_l\otimes {\mathcal O}_m, n=l+m,
\hat{S}_0=S\otimes 
{\mathcal O}_l$.
Recall that ${\mathcal O_n}={\mathbf F}[x_1,\ldots ,x_n]/(x_1^p,\ldots ,x_n^p)$ and $deg
x_i=-1$ with respect to the grading corresponding to the torus $T$. Therefore,
we can choose  the additional subalgebra ${\mathcal O}_l$ to be graded. In
such a case, ${\mathfrak n}_{-1}=
{\mathfrak l}_{-1}\oplus {\mathfrak  m}_{-1}$ where ${\mathfrak l}$ is the
maximal
ideal of ${\mathcal O}_l$.

According to (b) and (c), $S\otimes {\mathcal O}_n\subset {\mathfrak g}_0\subset 
(\Der S)\otimes {\mathcal O}_n + V + W_{n,0},$ where $V=1\otimes W_{n,1}.$ Recall that
${\mathcal O}_n={\mathbf F}[x_1,\ldots ,x_n]/(x_1^p,\ldots ,x_n^p), deg x_i=-1$ and 
$W_{n,1}=\, <\partial _1,\ldots ,\partial _n>\,$. This description is based on the
adjoint representation $\rho $ of  ${\mathfrak g}_0$ on the minimal ideal
$S\otimes {\mathcal O}_n$.
Note that $\rho $ is a faithful representation of ${\mathfrak g}_0$.
On the other hand,
$\hat{S}_0\otimes {\mathcal O}_m$ $\subset {\mathfrak g}_0$ $\subset 
\Der_0(\hat{S})\otimes {\mathcal O}_m + 1 $ $ \otimes $ $  W_m$.  To compare it with the
first description, consider the adjoint representation $\rho $ of ${\mathfrak
G}=\Der_0(\hat{S})\otimes {\mathcal O}_m + 1 $ $ \otimes $ $  W_m$ on $S\otimes {\mathcal
O}_n=\hat{S}_0\otimes {\mathcal O}_m$,

\begin{eqnarray} \rho ({\mathfrak G})
&\subset& \Der (\hat{S}_0)\otimes {\mathcal O}_m + 1\otimes
W_m \nonumber \\  &=& \Der (S\otimes {\mathcal O}_l)\otimes {\mathcal O}_m + 1\otimes W_m\nonumber \\ 
&=& (\Der S)\otimes {\mathcal O}_l\otimes {\mathcal O}_m + 1\otimes W_l\otimes 
{\mathcal O}_m + 1\otimes 1\otimes W_m.\nonumber\end{eqnarray}

\noindent Inasmuch as ${\mathfrak g}_0\subset \rho (\mathfrak G)$, we can, by taking into account the 
${\mathbb Z}$-grading of ${\mathfrak g}_0,$ obtain that

$${\mathfrak g}_0\subset (\Der S)\otimes {\mathcal O}_n + 1\otimes W_{n,1} +
1\otimes W_{l,1}\otimes {\mathfrak l}_{-1}+ 1\otimes W_{l,0}\otimes 1 + 1\otimes
1\otimes W_{m,0}.$$

\noindent According to (b), ${\mathcal O}_{n,-1}={\mathfrak n}_{-1}={\mathfrak
l}_{-1}\otimes 1\oplus
1\otimes{\mathfrak m}_{-1}$ is an irreducible $\pi (L_0)$-module where 
$\pi : (\Der S)\otimes 1 +1\otimes W_{n,0}\longrightarrow W_{n,0}$ is the
projection along the ideal $(\Der S)\otimes 1$. By (b), ${\mathfrak g}_0=S\otimes
{\mathcal O}_n +V+L_0,$ where $V=1\otimes W_{n,1}$. Hence,
 
$$\pi (L_0)\subset W_{l,1}\otimes {\mathfrak m}_{-1} +
W_{l,0}\otimes 1 + 1\otimes W_{m,0}.$$
\noindent This means that $1\otimes {\mathfrak m}_{-1}$ is invariant with respect to $\pi
(L_0)$. 
Thus, either ${\mathfrak m}_{-1}=0$ or ${\mathfrak m}_{-1}={\mathfrak n}_{-1};$
that is,
either $m=0$ or $m=n$.
 
Suppose that $m=0$. Then 

$$A(\hat{{\mathfrak g}})=\hat{S}=\hat{S}_{-u} + \ldots +\hat{S}_{-1} +\hat{S}_0
+\hat{S}_1+\ldots , \hat{S}_0=S\otimes {\mathcal O}_n.$$

Since $\hat{S}$ is simple, it follows that $\hat{S}_{-u}$ is an irreducible
$\hat{S}_0$-module. Furthermore,
as the nilradical $S\otimes {\mathfrak n}$ of $\hat{S}_0=S\otimes {\mathcal
O}_n$ has the 
${\mathbb Z}$-grading 
$S\otimes {\mathfrak n}=\oplus _{i>0}S\otimes {\mathcal O}_{n,-i}$ defined by
the torus $T$, it follows that it  acts as nilpotent elements on $\hat{{\mathfrak g}}$. Therefore, 
$[S\otimes {\mathfrak n},\hat{S}_{-u}]=0$. Since $\hat{S}_{-u}=\hat{{\mathfrak
g}}_{-u}$ is a $\hat{{\mathfrak g}}_0$-module, $\hat{{\mathfrak g}}_0={\mathfrak
g}_0$ and 
$\hat{S}_{0}=S\otimes {\mathcal O}_n$ is the minimal ideal of ${\mathfrak g}_0,$
it follows that $\hat{S}_{0}$ is contained in the kernel of the adjoint representation of
${\mathfrak g}_0$ on
$\hat{S}_{-u};$ that is, $[\hat{S}_{0},\hat{S}_{-u}]=0$. Since $\hat{S}$ is a
simple Lie algebra, 
it follows that $[\hat{S}_{-1},\hat{S}_1]=\hat{S}_0.$ Then

$$[\hat{S}_{-1},[\hat{S}_1,\hat{S}_{-u}]]=[[\hat{S}_{-1},\hat{S}_1],\hat{S}_{-u}
]=
[\hat{S}_0,\hat{S}_{-u}]=0.$$

\noindent Inasmuch as ${\mathfrak g}_{-1}=\hat{S}_{-1}$ and $\hat{{\mathfrak g}}$ is a
transitive Lie algebra, we have by Lemma 6 of [1] that if $x\in \hat{{\mathfrak
g}}_i, i>-u$ and $[\hat{{\mathfrak g}}_{-1},x]=0$, then $x=0$. Thus,
$[\hat{S}_1,\hat{S}_{-u}]=0.$
Set $\hat{S}^-\eqdef\oplus _{i\leqq 0}\hat{S}_i,$ and $\hat{S}^+\eqdef\oplus
_{i>0}\hat{S}_i$. Since $\hat{S}$ is simple, it follows that it is covered by the induced
module

$$ind \hat{S}_{-u}=U(\hat{S})\otimes_{U(\hat{S}^-)}\hat{S}_{-u}=
U(\hat{S}^+)\otimes _{\mathbf F}\hat{S}_{-u}.$$

\noindent Hence, $\hat{S}_{-u+1}=[\hat{S}_1,\hat{S}_{-u}]=0$. This contradiction shows that
$m=n$.
Thus, 

$$A(\hat{{\mathfrak g}})=\hat{S}\otimes {\mathcal O}_n, \hat{S}_0=S, 
\hat{{\mathfrak g}}_{-1}={\mathfrak g}_{-1}=\hat{S}_{-1}\otimes {\mathcal O}_n$$

\noindent and  $\hat{S}_1\otimes {\mathcal O}_n\subset \hat{{\mathfrak g}}_1={\mathfrak
g}_1\subset \Der_1(\hat{S})\otimes {\mathcal O}_n$.
Therefore, the bracket operation in ${\mathfrak g}$, 
$[ , ]:{\mathfrak g}_{-1}\times {\mathfrak g}_1\longrightarrow [{\mathfrak
g}_{-1},{\mathfrak g}_1]$ is the restriction of the ${\mathcal O}_n$-bilinear
bracket operation in 
$\Der(\hat{S})\otimes {\mathcal O}_n$,

$$\hat{S}_{-1}\otimes {\mathcal O}_n\times \Der_1(\hat{S})\otimes {\mathcal O}_n
\longrightarrow \hat{S}_0\otimes {\mathcal O}_n=S\otimes {\mathcal O}_n=
[{\mathfrak g}_{-1},{\mathfrak g}_1].$$

\bigskip

\noindent  {\it (i) Conclusion of the proof of the Theorem\/.}

Let ${\mathfrak g}_1=\oplus _{i\geqq k}{\mathfrak g}_{1,-i}$. Suppose that
$k\leqq |\delta |=n(p-1)$. By (h) the bracket operation $[{\mathfrak g}_1,
{\mathfrak g}_{-1}]$ is 
a restriction of the ${\mathcal O}_n$-bilinear bracket operation. Therefore,
$0\neq [{\mathfrak g}_{1,-k}, {\mathfrak g}_{-1}]=
[{\mathfrak g}_{1,-k}, \hat{S}_{-1}\otimes {\mathcal O}_n]$ is  an ${\mathcal
O}_n$-submodule of $[{\mathfrak g}_1,{\mathfrak g}_{-1}]=S\otimes {\mathcal
O}_n$.
On the other hand, by (d),

\begin{eqnarray} [{\mathfrak g}_{1,-k},{\mathfrak g}_{-1}] &=& [{\mathfrak g}_{1,-k}, {\mathfrak
g}_{-1,1}+
\ldots +{\mathfrak g}_{-1,|\delta |+1}]\nonumber \\ 
&\subset& {\mathfrak g}_{0,-k+1}+{\mathfrak g}_{0,-k+2}+\ldots\nonumber \\
&=& S\otimes {\mathcal O}_{n,-k+1}+\ldots +S\otimes 1.\nonumber\end{eqnarray}

\noindent Hence, 
$S\otimes {\mathcal O}_{n,-|\delta |}\cap [{\mathfrak g}_{1,-k}, {\mathfrak
g}_{-1}]=0$.
However, any nonzero ${\mathcal O}_n$-submodule of $S\otimes {\mathcal O}_n$ has
a nontrivial intersection with $S\otimes x^{\delta }=S\otimes {\mathcal O}_{n,
-|\delta |}$.
We have arrived at a contradiction. Consequently, $k=|\delta |+1$.  By (g), ${\mathfrak g}_{1,-|\delta
|-1}\cong 
L_{-|\delta |-1}$ and $[S, L_{-|\delta |-1}]\neq 0$.

We now show that if $U\subset L_{-|\delta |-1}$ and $[S, U]\neq 0$, then $[S\otimes x,
U]\neq 0$ in $L$. Here $x\in {\mathcal O}_{n,-1},$ and $S\otimes x\subset S\otimes
{\mathcal O}_{n,-1}=
L_{-1}$. 

Let $1\otimes \partial _{\xi }\in 1 \,\,\, \otimes <\partial _1, \ldots ,\partial
_n> \,\,\, = V\subset L_1, \xi (x)\neq 0$. By (f), 2), $M({\mathcal L}_0)=L_{-|\delta
|-1}+ L_{-|\delta |-2} +\ldots $. Hence
$[V, L_{-|\delta |-1}]=0$ and 

$$[1\otimes \partial _{\xi },[S\otimes x, U]]=[[1\otimes \partial _{\xi },
S\otimes x], U]=
[S, U]\neq 0.$$

\noindent Therefore, $[S\otimes x, U]\neq 0$ for any $x\in {\mathcal O}_{n,-1}$. Since
$[S, L_{-|\delta |-1}]\neq 0$ and 
$L_{-|\delta |-1}=[L_{-1}, L_{-|\delta |}]=\sum _{i=1}^n[S\otimes x_i,
L_{-|\delta |}],$
there exists an $i$ such that $[S, [S\otimes x_i, L_{-|\delta |}]]\neq 0$. Set 
$U=[S\otimes x_i, L_{-|\delta |}]\subset L_{-|\delta |-1}$. Then $[S, U]\neq
0.$ Thus, as above,
$[S\otimes x_j,U]=[S\otimes x_j,[S\otimes x_i, L_{-|\delta |}]]\neq 0$.
We may renumber the variables and  suppose $i=1$.  Since by  assumption $\dim
V=n>1,$ it follows that
$[S\otimes x_2,[S\otimes x_1, L_{-|\delta |}]]\neq 0$ in $L$.
For $1\otimes \partial _2\in V\subset L_1,$ we have

\begin{eqnarray} [[S\otimes x_2,[S\otimes x_1, L_{-|\delta |}]], 1\otimes \partial _2] &=& 
[[S\otimes x_2, 1\otimes \partial _2], [S\otimes x_1, L_{-|\delta |}]]\nonumber \\
&=& [S,[S\otimes x_1, L_{-|\delta |}]]\neq 0.\nonumber\end{eqnarray}

\noindent On the other hand, $[S\otimes x_1,1\otimes \partial _2]=0$ and 

$$[[S\otimes x_2, L_{-|\delta |}], 1\otimes \partial _2]\subset [L_{-|\delta
|-1},V]=0.$$

\noindent Inasmuch as $M({\mathcal L}_0)=L_{-|\delta |-1}+L_{-|\delta |-2}+\ldots $, we have that
$L_{-i}=S\otimes {\mathcal O}_{n,-i}, i=1,\ldots , |\delta |$ and the
bracket operation in $L$, 
$L_{-i}\times L_{-j}\longrightarrow L_{-i-j}$ coincides with the bracket operation
in 
$S\otimes {\mathcal O}_n$ for any $i, j$ such that $1\leqq i, j, i+j \leqq |\delta
|$.  Hence, in $L,$

\begin{eqnarray} 
[S\otimes x_1x_2, S\otimes x^{\delta }] &=& [S\otimes x_1x_2, [S\otimes x_1^{p-1},S\otimes x_2^{p-1}\ldots
x_n^{p-1}]]\nonumber \\
&=& [[S\otimes x_1x_2, S\otimes x_1^{p-1}], S\otimes x_2^{p-1}\ldots x_n^{p-1}]\nonumber\\ 
&&\qquad + \, [S\otimes x_1^{p-1}, [S\otimes x_1x_2, 
S\otimes x_2^{p-1}\ldots x_n^{p-1}]]\nonumber \\
&=& 0.\nonumber\end{eqnarray}

\noindent Now, as $[S\otimes x_1,1\otimes \partial _2]=0$ and 
$[[S\otimes x_2, L_{-|\delta |}], 1\otimes \partial _2]\subset [L_{-|\delta
|-1},V]=0$, we have

\begin{eqnarray}
[[S\otimes x_2, [S\otimes x_1, L_{-|\delta |}]], 1\otimes \partial_2] 
&=& [[[S\otimes x_2, S\otimes x_1], L_{-|\delta |}], 1\otimes \partial _2]\nonumber\\
&&\qquad + \, [[S\otimes x_1,[S\otimes x_2, L_{-|\delta |}]], 1\otimes \partial _2]\nonumber\\
&=&
[[S\otimes x_1x_2, L_{-|\delta |}], 1\otimes \partial _2]\nonumber \\
&=& [[S\otimes x_1x_2, S\otimes x^{\delta }], 1\otimes\partial _2]=0.\nonumber
\end{eqnarray}

\noindent The contradiction obtained completes the proof of the theorem. \qed

\begin{center} References
\end{center}

\smallskip

\noindent [1] {\it Benkart G.M., Gregory T.B.}  Graded Lie
algebras with classical reductive null component,  Math. Ann., 1989, 285, pp. 85--98.

\smallskip

\noindent [2] {\it  Benkart G.M.,  Gregory T.B., 
Premet A.A.}  The Recognition Theorem for Graded Lie Algebras in Prime
Characteristic, Memoirs of the American Mathematical Society, Volume 197, Number 120, 2009.

\smallskip

\noindent [3] {\it Benkart G.M., Kostrikin A.I., Kuznetsov M.I.}  The simple
 graded Lie algebras of characteristic three with classical reductive component 
$L_0$, Comm. in Algebra, 1996, 24, pp. 223--234.

\smallskip

\noindent [4] {\it Block R.E.} Modules over differential polynomial rings, Bull. Amer. Math. Soc., 1973, 79, pp. 729 -- 733.

\smallskip

\noindent [5] {\it Brown G.}  On the structure of some Lie algebras of Kuznetsov, Michigan Math. J., 1992, 39, pp. 85--90.

\smallskip

\noindent [6]  {\it Curtis C.W.}  Representations of Lie 
algebras of classical type with applications to linear groups,
J. Math. Mech., 1960, 9, pp. 307--326.

\smallskip

\noindent [7]  {\it  Gregory T.B.,  Kuznetsov M.I.}  On
depth-three graded Lie algebras of characteristic three with
classical reductive null component, Comm. in Algebra,
2004, 33(9), pp. 3339--3371.

\smallskip

\noindent [8]  {\it  Kac V.G.} The classification of simple Lie algebras over a field of nonzero characteristic,
Izv. Akad. Nauk SSSR, Ser. Mat., 1970, 34(2), pp. 385-408 (Russian); English transl., Math. USSR-Izv, 1970, 4, pp. 39--413.
  
\smallskip

\noindent [9] {\it  Kostrikin A.I., Ostrik V.V.} To the recognition theorem for Lie algebras of
characteristic 3, Mat. Sbornik, 1995, 186(5),  pp. 73--88 (Russian).

\smallskip

\noindent [10] {\it  Kostrikin A.I., Shafarevich I.R.} Graded Lie algebras of finite characteristic, Izv. Akad. Nauk
SSSR, Ser. Mat., 1969, 33(2), pp. 251--322 (Russian); English transl., Math. USSR-Izv., 1969, 3, pp. 237--304.

\smallskip

\noindent [11] {\it Kuznetsov M.I.} Truncated induced modules over transitive Lie algebras  of characteristic p, Izv. Akad. Nauk. SSSR, Ser. Mat., 1989, 53(3), pp. 557--589 (Russian); 
English transl., Math. USSR-Izv., 1990, 34, pp. 575--608.

\smallskip

\noindent [12] {\it Skryabin S.M.} New series of simple Lie algebras of characteristic 3,Mat. Sbornik, 
1992, 183, pp. 3-22 (Russian); English transl., Russian Acad. Sci. Sb. Math., 1993, 70, pp. 389--406.

\smallskip

\noindent [13]  {\it  Strade H.}  Simple Lie algebras over Fields of Positive Characteristic: I. Structure
Theory, DeGruyter Expositions in Mathematics V. 38, New York: Walter de
Greyter, 2004.

\smallskip

\noindent [14] {\it Weisfeiler B.} On filtered Lie algebras and their associated graded algebras,
Functional Anal. i. Prilozhen., 1968, 2(1), pp. 94--95 (Russian).

\smallskip

\noindent [15] {\it Weisfeiler B.}  On the structure of the
minimal ideal of some graded Lie algebras in characteristic $p >0$, J. Algebra, 1978, 
53(2), pp. 344--361.

\smallskip

Department of Mathematics, The Ohio State University at Mansfield,
Mansfield, Ohio  44906, USA

\smallskip

Nizhny Novgorod State University, Gagarin Ave., 23
Nizhny Novgorod, 603095, Russia

\end{document}